\documentclass[12pt,a4paper]{article}
\pdfoutput=1
\usepackage[english]{babel}
\usepackage{verbatim}
\usepackage{bm}

\usepackage{amsmath}
\usepackage{amssymb}
\usepackage{amsfonts}
\usepackage{amsbsy}

\usepackage{graphicx}
\usepackage{graphics}
\usepackage{subfigure}
\usepackage{epsfig}

\usepackage{array}
\usepackage{booktabs}% for top-,mid- & bottomrule
\usepackage{caption}%for \captionof
\usepackage{tabularx}

\renewcommand{\P}{\mathbb P}
\newcommand{\Z}{\mathbb Z}

\begin{document}

\title{Flux form Semi-Lagrangian methods \\
for parabolic problems}

\author{Luca Bonaventura$^{(1)}$,\ \  Roberto Ferretti$^{(2)}$}

\maketitle

\begin{center}
{\small
$^{(1)}$ MOX -- Modelling and Scientific Computing, \\
Dipartimento di Matematica, Politecnico di Milano \\
Via Bonardi 9, 20133 Milano, Italy\\
{\tt luca.bonaventura@polimi.it}
}
\end{center}

\begin{center}
{\small
$^{(2)}$
Dipartimento di Matematica e Fisica\\
Universit\`a degli Studi Roma Tre\\
 L.go S. Leonardo Murialdo 1,
00146, Roma, Italy\\
{\tt ferretti@mat.uniroma3.it}
}
\end{center}

\date{}


\vspace*{0.5cm}


\vspace*{0.5cm}

\pagebreak

\abstract{A semi-Lagrangian method for parabolic problems is proposed,
that extends previous work by the authors to achieve a fully conservative, flux-form 
discretization of linear and nonlinear diffusion equations. 
A basic consistency and convergence analysis are proposed. Numerical examples
validate the proposed method and display its potential for consistent semi-Lagrangian discretization
of advection--diffusion and nonlinear parabolic problems.} 

\pagebreak

 \section{Introduction.}
 \label{intro} \indent
 
 It is well known that, in their original formulation, Semi-Lagrangian (SL) schemes were developed
 for the case of purely hyperbolic equations and do not guarantee conservation of mass.
  In recent years, an increasing number of
SL discretization approaches have been proposed to circumvent these limitations.

Various extensions of the SL strategy to second order, parabolic problems, appear for example in~\cite{camilli:1995},~\cite{ferretti:2010},~\cite{milstein:2002},~\cite{milstein:2000},~\cite{milstein:2001},~\cite{teixeira:1999} (see also~\cite{falcone:2013} for a more complete review),
as opposite to the ``classical'' technique of treating in Eulerian form diffusion terms within SL schemes.
In spite of the variety of applications involved, the common feature of these works is to replace the  computation of the solution at the foot of a characteristic with a weighted average of values computed at multiple points.
A rigorous derivation of these techniques can be based on the Feynman--Ka\v{c} representation formula, with stochastic trajectories integrated backward in time as in conventional SL schemes. However, due to their stochastic origin, all these generalizations of the SL approach are based on the trace form of parabolic problems, and unsuitable for the divergence form
\begin{equation}
u_t  =  \text{div}\left (D\nabla u \right) \label{divform}
\end{equation}
which is the most widely used in many applications, especially those to computational fluid dynamics on geophysical scales, for which SL methods for advective problems have proven to be especially useful  (see e.g. the review in~\cite{bonaventura:2012}). In the previous paper~\cite{bonaventura:2014}, we have proposed for the first time a SL method for parabolic problems in divergence form. Numerical experiments reported therein show that the method, albeit formally only first-order accurate in time, allows to compute remarkably accurate approximations of linear and nonlinear parabolic problems, as well as to achieve easily higher order accuracy in space also on nonuniform meshes. 

 On the other hand, the last two decades have witnessed an increasing interest in conservative, flux-form SL (FFSL) schemes which have been proposed, e.g., in  ~\cite{crouseilles:2010},~\cite{frolkovic:2002},~\cite{leonard:1996},~\cite{lin:1996},~\cite{phillips:2001},
 ~\cite{qiu:2011a},~\cite{qiu:2011b},~\cite{restelli:2006},~\cite{rossmanith:2011},
 for various linear and nonlinear advective problems. These methods are
 exactly mass conservative, as opposed to methods based on upstream remapping of computational cells
 (also called inherently conservative in the literature), see e.g.    ~\cite{dukowicz:1984},~\cite{dukowicz:2000}, 
 ~\cite{laprise:1995},~\cite{leslie:1995},~\cite{lipscomb:2005}, ~\cite{machenhauer:1997},~\cite{nair:2002a},~\cite{rancic:1995},
 ~\cite{zerroukat:2002},~\cite{zerroukat:2004},~\cite{zerroukat:2005},~\cite{zerroukat:2009}.

The method proposed in~\cite{bonaventura:2014}, however, even though discretizing Equation~\eqref{divform} directly, is not formulated in flux form at the discrete level and  does not guarantee exact mass conservation. 
In the present work, we try to fill the gap by extending the approach of~\cite{bonaventura:2014} to achieve a fully conservative, flux form discretization of the diffusion equation~\eqref{divform}. Along the same lines as~\cite{bonaventura:2014}, we outline a consistency and stability analysis of the method in a simplified setting, and perform a number of numerical simulations to validate the proposed method.

The paper is structured as follows. In Section~\ref{second_ord}, our novel FFSL discretization is introduced, while an analysis of its consistency and stability properties will be presented in Sections~\ref{consistency} and~\ref{stability}. Section~\ref{ndim} treats the extensions of the scheme to problems in higher dimensions. Finally, numerical results obtained with the proposed method for both linear and nonlinear models will be presented in Section~\ref{tests}, while some conclusions on the potential advantages of our approach will be drawn in Section~\ref{conclu}.

\section{Flux Form Semi-Lagrangian methods for second order problems.}
\label{second_ord}
In order to sketch the basic ideas of the scheme, we restrict for the moment to the approximation of the pure diffusion equation
\begin{equation}\label{diffusione}
u_t = (\nu(x,t) u_x)_x
\end{equation}
in a single space dimension, considered for simplicity on the whole real line. 
The extension to advection--diffusion equations and to multiple dimensions will be discussed later in this Section and in Section~\ref{ndim}. We have shown in~\cite{bonaventura:2014} that~\eqref{diffusione} can be approximated via an abstract difference operator in the form of a convex combination of pointwise values. More specifically,
let $\Delta t$ and $\Delta x$ denote the time and space discretization steps, respectively, with $t_n=n\Delta t$ for $n\in [0,T/\Delta t]$. The space grid is supposed to be infinite and uniform, that is, for $i\in\Z,$ $x_i= i\Delta x$.
In the case of the conservative scheme, we will also consider the intermediate points
$x_{i\pm 1/2}= (i\pm 1/2)\Delta x, $
which will appear as endpoints of grid cells $\Omega_i=[x_{i-1/2},x_{i+1/2}]$.  The nonconservative method is derived by setting
\begin{equation}\label{schema}
u(x_i,t_{n+1}) \approx A_i^+ u(x_i+\delta_i^+,t_n) + A_i^- u(x_i-\delta_i^-,t_n),
\end{equation}
provided the constants $A_i^+$, $A_i^-$, $\delta_i^+$, $\delta_i^+$ satisfy the consistency conditions
\begin{equation}\label{sistema}
\begin{cases}
A_i^+ + A_i^- = 1 + O(\Delta t^2) \\
A_i^+\delta_i^+ - A_i^-\delta_i^- = \Delta t \> \nu_x + O(\Delta t^2) \\
A_i^+{\delta_i^+}^2 + A_i^-{\delta_i^-}^2 = 2\Delta t \> \nu + O(\Delta t^2) \\
A_i^+{\delta_i^+}^3 - A_i^-{\delta_i^-}^3 = O(\Delta t^2).
\end{cases}
\end{equation}
Note that, actually, the additional condition ${\delta_i^\pm}^4 = O(\Delta t^2)$ should also be added, which is however
 already implied by the second condition. We will denote the numerical solutions  at time $t_n$ by the vector $V^n=(v_i^n)_{i\in\Z}$. While in the nonconservative form SL scheme these values should be understood as pointwise values, in a FFSL scheme they must be interpreted as cell averages. 

In~\cite{bonaventura:2014}, the abstract structure introduced in~\eqref{sistema} takes the specific form of the nonconservative approximation
\begin{equation}\label{schema_tm}
v_i^{n+1} = \frac{1}{2}I_p[V^n](x_i+\delta^+_i) + \frac{1}{2}I_p[V^n](x_i-\delta^-_i),
\end{equation}
where the pointwise values are reconstructed by the interpolation operator (of degree $p$) $I_p[V^n]$,
while the displacements $\delta_i^\pm$ are given as solution of the equation
\begin{equation}\label{delta_t}
\delta^\pm_i = \sqrt{2\Delta t\>\nu(x_i\pm\delta^\pm_i,t^n)}.
\end{equation}
For a flux form variant of the above approach, rather than approximating a pointwise value, the cell averages
\[
\bar u_i(t) \approx \frac{1}{\Delta x} \int_{\Omega_i} u(x,t) dx
\]
of the solutions of Equation~\eqref{diffusione} must be approximated, as customary in finite volume
methods, to obtain the  discrete mass balance
\begin{equation}\label{ff}
\bar u_i(t_{n+1}) \approx \bar u_i(t_n) + \frac{1}{\Delta x}\left(\mathcal{F}_{i+1/2}^n-\mathcal{F}_{i-1/2}^n\right).
\end{equation}
The flux $\mathcal{F}_k^n$ might be defined in abstract form as
\[
\mathcal{F}_k^n = \frac{1}{2}\left( \int_{x_k}^{x_k+\delta_k}u(x,t_n)dx -  \int_{x_k-\delta_k}^{x_k}u(x,t_n)dx  \right)
\]
for a properly defined $\delta_k$. The practical version of the scheme implements this operator with a numerically reconstructed function, 
so that the proposed FFSL method for  Equation~\eqref{diffusione} is defined as
\begin{equation}\label{schema_ff}
v_i^{n+1} = v_i^{n}+\frac{1}{\Delta x}\left(F_{i+1/2}^n-F_{i-1/2}^n\right)
\end{equation}
with a numerical flux given, for $k=i\pm 1/2$, by
\begin{equation}\label{flux}
F_k^n = \frac{1}{2}\left(  \int_{x_k}^{x_k+\delta_k}R_q[V^n](x)dx- \int_{x_k-\delta_k}^{x_k}R_q[V^n](x)dx  \right),
\end{equation}
and with
\begin{equation}\label{delta_ff}
\delta_k = \sqrt{2\Delta t\>\nu(x_k,t^n)}.
\end{equation}
Here, the operator $R_q[V]$ is a polynomial reconstruction of degree $q$ which satisfies over each cell $\Omega_m=[x_{m-1/2},x_{m+1/2}]$ the properties of having the correct cell average $w_m$ and reconstructing with order $O(\Delta x^{q+1})$ a smooth solution. There is a close relationship between the two operators $I_p$ and $R_q$. We refer to~\cite{ferretti:2013b} for an in-depth discussion of these theoretical issues, while in Section~\ref{stability} will recall the main points of interest for the purpose of a stability analysis of the proposed method.
It should be remarked that, in this case, it is not required to solve any equation to obtain $\delta_k$.
Furthermore, notice that both methods~\eqref{schema_tm} and~\eqref{schema_ff} 
employ viscosity values frozen at the time level $t_n$, so that an extension to the nonlinear case
is straightforward.
 
%At a first glance, it is not obvious that~\eqref{flux} retains the correct physical meaning of a flux. 
In order to provide a more intuitive interpretation of the proposed scheme,
 it could be observed that the effect of diffusion (as modeled by Fick's law implicit in  
Equation~\eqref{diffusione})
 is to move mass from regions with high density towards regions with low density. 
 The difference of integrals in~\eqref{flux} can be shown to be an approximate model of this process.
 Indeed, the interval $\sqrt{2\Delta t\>\nu}$ is   the correct space scale associated to the given diffusion over a time step $\Delta t$.
Furthermore, it can be observed that
\[
\int_{x_k}^{x_k+\delta_k}u(x,t_n)dx \approx  \delta_k\> u\left(x_k+\frac{\delta_k}{2},t_n\right),
\]
\[
\int_{x_k-\delta_k}^{x_k}u(x,t_n)dx \approx  \delta_k\> u\left(x_k-\frac{\delta_k}{2},t_n\right),
\]
so that, as a consequence, one obtains
\begin{eqnarray*}
\mathcal{F}_k^n & \approx & \frac{\delta_k}{2}\left( u\left(x_k+\frac{\delta_k}{2},t_n\right) - u\left(x_k-\frac{\delta_k}{2},t_n\right) \right) \approx \\
& \approx & \frac{\delta_k^2}{2}\> u_x(x_k,t_n) =  \Delta t\>\nu(x_k,t^n)\> u_x(x_k,t_n)
\end{eqnarray*}
which clearly makes~\eqref{ff} consistent with~\eqref{diffusione}. Of course, this is rather a heuristic explanation than a rigorous consistency analysis. Consistency (and its optimal rate) will be proven in the next Section on the basis of conditions~\eqref{sistema}.

We now show briefly how the proposed FFSL method can be combined with analogous methods
for the flux form advection equation. Consider the advection-diffusion equation
\begin{equation}
\label{advdiff}
u_t + (f(x,t)u)_x = (\nu(x,t) u_{x} )_x,
\end{equation}
taken again for simplicity on the infinite real line. The combination of the
FFSL methods for advection and diffusion is obtained, for the purposes of this paper, by a simple operator splitting
approach. This combination results in a method that is first order accurate in time,  which is compatible with the
time accuracy of the FFSL method for diffusion, as it will be seen in Section~\ref{consistency}.
 In a first step, any FFSL method for advection can be used. In Section~\ref{tests}, we consider for concreteness the well-known method  presented in~\cite{lin:1996}. This results in a numerical approximation of the cell averages ${\tilde v}_i^{n+1}$ of the solution of~\eqref{advdiff} with $\nu=0.$ An approximation of the
solution of the complete Equation~\eqref{advdiff}  is then obtained by computation of
formula~\eqref{schema_ff} based on the intermediate values ${\tilde v}_i^{n+1}.$
 
% \pagebreak

\section{Consistency.}
\label{consistency}

To prove consistency of the flux form SL method, we start by neglecting the error associated to the 
space discretization and rewriting the right-hand side of~\eqref{ff} in the more convenient form:
\begin{eqnarray*}
\bar u_i(t_n) + \frac{1}{\Delta x}\left(\mathcal{F}_{i+1/2}^n-\mathcal{F}_{i-1/2}^n\right) & = &
 \bar u_i(t_n) + \frac{1}{2\Delta x}\left(\int_{x_{i+1/2}}^{x_{i+1/2}+\delta_{i+1/2}}- \int_{x_{i+1/2}-\delta_{i+1/2}}^{x_{i+1/2}}  - \right.\nonumber\\
&& \left. - \int_{x_{i-1/2}}^{x_{i-1/2}+\delta_{i-1/2}} + \int_{x_{i-1/2}-\delta_{i-1/2}}^{x_{i-1/2}} \right)  u(x,t_n) dx = \\
& = & \bar u_i(t_n) + \frac{1}{2\Delta x}\left(\int_{x_{i-1/2}-\delta_{i-1/2}}^{x_{i+1/2}-\delta_{i+1/2}} - \int_{x_{i-1/2}}^{x_{i+1/2}} + \right. \nonumber \\
&& \left. + \int_{x_{i-1/2}+\delta_{i-1/2}}^{x_{i+1/2}+\delta_{i+1/2}} - \int_{x_{i-1/2}}^{x_{i+1/2}} \right) u(x,t_n) dx. \nonumber 
\end{eqnarray*}
We obtain therefore
\begin{equation}\label{ff2}
\bar u_i(t_n) + \frac{1}{\Delta x}\left(\mathcal{F}_{i+1/2}^n-\mathcal{F}_{i-1/2}^n\right) = \frac{1}{2\Delta x}\left(\int_{x_{i-1/2}-\delta_{i-1/2}}^{x_{i+1/2}-\delta_{i+1/2}} + \int_{x_{i-1/2}+\delta_{i-1/2}}^{x_{i+1/2}+\delta_{i+1/2}} \right) u(x,t_n) dx.
\end{equation}
The terms $\delta_{i\pm 1/2}$ can also be expressed as
\begin{eqnarray}\label{delta_ff2}
\delta_{i\pm 1/2} & = & \sqrt{2\Delta t \> \nu(x_{i\pm 1/2})} = \\
& = & \sqrt{2\Delta t \left( \nu \pm \frac{\Delta x}{2}\nu_x + \frac{\Delta x^2}{8}\nu_{xx} + O\left(\Delta x^3\right) \right)} =  \nonumber \\
& = & \sqrt{2\Delta t \> \nu \left( 1 \pm \frac{\Delta x\>\nu_x}{2\nu} + \frac{\Delta x^2\>\nu_{xx}}{8\nu} + O\left(\Delta x^3\right) \right)} = \nonumber \\
& = & \delta_i \left( 1 \pm \frac{\Delta x\>\nu_x}{4\nu} + \frac{\Delta x^2}{16\nu}\left(\nu_{xx} - \frac{\nu_x^2}{2\nu^2} \right) + O\left(\Delta x^3\right)\right),\nonumber
\end{eqnarray}
where, in the last row, we have applied the Taylor expansion of the square root and collected all terms of order $O(\Delta x^3)$.
It is now easy to recognize that, once replaced (up to a $O(\Delta x^2)$ error) the integrals in the last row of~\eqref{ff2} with their midpoint approximation, i.e.,
\[
\int_{x_{i-1/2}\pm\delta_{i-1/2}}^{x_{i+1/2}\pm\delta_{i+1/2}} u(x,t_n) dx = (\Delta x \pm \delta_{i+1/2} \mp \delta_{i-1/2}) u(x_i\pm\delta_i) + O\left(\Delta x^2\right),
\]
and using~\eqref{delta_ff2}, then~\eqref{ff2} is in the form~\eqref{schema}, with
\begin{eqnarray}\label{A_delta_ff}
\delta_i^{\pm} & = & \delta_i \label{A_delta_ff1} \\
A_i^\pm & = & \frac{\Delta x \pm (\delta_{i+1/2} - \delta_{i-1/2})}{2\Delta x} =  \\
& = & \frac{1}{2}\pm \frac{\delta_i \nu_x}{4\nu} + O\left(\Delta x^2\right).\label{A_delta_ff2}\nonumber
\end{eqnarray}
We check now that~\eqref{A_delta_ff1}--~\eqref{A_delta_ff2} satisfy the consistency conditions~\eqref{sistema}. The first condition is clearly satisfied. As for the second, we have
\begin{eqnarray*}
A_i^+\delta_i^+ - A_i^-\delta_i^- & = & \delta_i\left(A_i^+ - A_i^- \right) = \\
& = & \frac{\delta_i^2\nu_x}{2\nu} + O\left(\Delta x^2\right) = \\
& = & \Delta t\>\nu_x + O\left(\Delta x^2\right).
\end{eqnarray*}
The third condition is rewritten as
\begin{eqnarray*}
A_i^+{\delta_i^+}^2 + A_i^-{\delta_i^-}^2 & = & \delta_i^2\left(A_i^+ + A_i^- \right) = \\
& = & 2\Delta t\>\nu\left(1 + O\left(\Delta x^2\right)\right).
\end{eqnarray*}
Finally, for the fourth condition we obtain
\begin{eqnarray*}
A_i^+{\delta_i^+}^3 - A_i^-{\delta_i^-}^3 & = & (2\Delta t\>\nu)^{3/2} \left(\frac{\delta_i\nu_x}{2\nu} + O\left(\Delta x^2\right) \right) = \\
& = & O\left(\Delta t^2\right),
\end{eqnarray*}
in which the last expression is  achieved by taking into account that
$$
\frac{\delta_i\nu_x}{2\nu} =O\left(\Delta t^{1/2}\right).
$$
Therefore, the abstract operator is consistent.
Finally, introducing the polynomial reconstruction $R_q[V]$, we note that fluxes are approximated with accuracy of order $O(\Delta x^{q+2})$. Then, the estimate 
\begin{equation}\label{cons_err}
L(\Delta x, \Delta t) = O\left(\Delta t^r + \frac{\Delta x^s}{\Delta t} \right),
\end{equation}
holds for the truncation error $L,$ with $r=1$ and $s=q+2$.

\section{Stability.}
\label{stability}

We briefly discuss the stability of the  proposed scheme in the case of a constant coefficient equations,
 while variable coefficient equations seem to require a deeper study.
First, we note that a well understood case is the construction of the nonconservative scheme with $I_p$ in the form of a symmetric Lagrange interpolation, for $p$ odd. For example, if $p=3$, the value $I_3[V] (x)$ is computed, for $x\in [x_k,x_{k+1}]$, by interpolating the values $v_{k-1},\ldots,v_{k+2}$. In the pure advection case, this construction is known to be stable (see~\cite{ferretti:2010b}).
In order to construct a stable conservative scheme, we can define the reconstruction $R_q$ (see~\cite{ferretti:2013b} for more details and examples) according to the axioms

\begin{itemize}

\item[i)]
For $x\in \Omega_m$, $R_q[W](x)=Q_m(x)$ for some polynomial $Q_m\in \P_q$, with $q$ even;

\item[ii)]
For $x\in \Omega_m$, $Q_m(x)$ depends on the values $w_k$, with $k=m-q/2,\ldots,m+q/2$, and more precisely it satisfies the conditions
\begin{equation*}
\frac{1}{\Delta x} \int_{\Omega_k} Q_m(x)dx = w_k \qquad (k=m-q/2,\ldots,m+q/2).
\end{equation*}

\end{itemize}
With this definition, it is proven in~\cite{ferretti:2013b} that, if $p=q+1$, then
\begin{equation}\label{moving_av}
\frac{1}{\Delta x}\int_{x-\Delta x/2}^{x+\Delta x/2} R_q[W](\xi)d\xi = I_p[W](x).
\end{equation}

Consider now the problem~\eqref{diffusione} for a constant diffusivity $\nu$. In this case $\delta_i\equiv\delta$, so that, rewriting~\eqref{ff2} for the reconstructed numerical solution and using~\eqref{moving_av}, we get
\begin{eqnarray}\label{ff3}
v_i^{n+1} & = &v_i^n + \frac{1}{\Delta x}\left(F_{i+1/2}^n-F_{i-1/2}^n\right) =
 \\
& = & \frac{1}{2}\left(\frac{1}{\Delta x}\int_{x_{i-1/2}-\delta}^{x_{i+1/2}-\delta} + \frac{1}{\Delta x}\int_{x_{i-1/2}+\delta}^{x_{i+1/2}+\delta} \right) R_q[V^n](x) dx = \nonumber \\
& = & \frac{1}{2}I_p[V^n](x-\delta) + \frac{1}{2}I_p[V^n](x+\delta).\nonumber
\end{eqnarray}
The scheme can then be recast as the convex combination
\begin{equation}\label{media}
V^{n+1} = \frac{1}{2} B^+ V^n + \frac{1}{2} B^- V^n,
\end{equation}
where the terms $B^\pm V^n$ represent respectively the schemes written in extended form as
\begin{equation*}
v_i^{n+1}=I_p[V^n] (x_i\pm\delta),
\end{equation*}
which (see~\cite{ferretti:2010b}) are stable in the 2-norm for any $\delta$ and $p$. This implies that
\[
\left\| B^\pm \right\|_2 \le 1,
\]
and therefore that the complete scheme is also stable in the same norm.

% \section{Coupling to flux form SL methods for advection.}
%\label{couple}

 \section{Multiple space dimensions.}
\label{ndim}
In this Section, we discuss the extension of the proposed approach to the $d$-dimensional case. 
The extension is straightforward in the case of a structured orthogonal grid and diagonal diffusivity matrix
\begin{equation}\label{diff_diag}
\Lambda(x,t) = \text{diag} (\nu_1(x,t),\ldots,\nu_d(x,t)),
\end{equation}
and in particular, for a variable but isotropic diffusion (for which $\nu_1(x,t)=\cdots =\nu_d(x,t)$). Then, the diffusion equation reads
\begin{eqnarray*}
u_t & = & \text{div}(\Lambda(x,t)\nabla u) = \\
& = & \sum_{j=1}^d (\nu_j(x,t)u_{x_j})_{x_j}.
\end{eqnarray*}
By a first-order expansion of the solution with respect to time, we get
\begin{eqnarray*}
u(x,t+\Delta t) & = & u(x,t) + \Delta t \sum_{j=1}^d (\nu_j u_{x_j})_{x_j}(x,t) + O(\Delta t^2) = \\
& = & \sum_{j=1}^d \left[ \frac{u(x,t)}{d}+\frac{\Delta t}{d} (d\nu_j u_{x_j})_{x_j}(x,t)\right] + O(\Delta t^2) = \\
& = &  u(x,t) + \frac{1}{d} \sum_{j=1}^d \left[\Delta t (d\nu_j u_{x_j})_{x_j}(x,t)\right] + O(\Delta t^2),
\end{eqnarray*}
which shows that, up to first-order accuracy, the $d$-dimensional version can be obtained by averaging the diffusion operators in each direction, with the only modification of scaling each one-dimensional diffusivity by a factor $d$. This allows to split the diffusion along each of the variables, where (in case of an orthogonal mesh) flux tubes are also aligned with the grid.

For example, the 2-dimensional equation
\[
u_t = \left(\nu_1(x,t)u_{x_1}\right)_{x_1} + \left(\nu_2(x,t)u_{x_2}\right)_{x_2}
\]
would be approximated by the scheme (written with an obvious notation at the node $x_i=(x_{i_1},x_{i_2})$)
\begin{equation}\label{eq:schema2d}
v_i^{n+1} = v_i^n + \frac{1}{\Delta x^2}\left(F_{i_1+1/2,i_2}^n-F_{i_1-1/2,i_2}^n+F_{i_1,i_2+1/2}^n-F_{i_1,i_2+1/2}^n\right),
\end{equation}
in which, for example, the flux $F_{i_1+1/2,i_2}^n$ would be given by
\begin{eqnarray}\label{flux2d}
F_{i_1+1/2,i_2}^n & = & \frac{1}{4}\left(\int_{x_{i_2-1/2}}^{x_{i_2+1/2}}\int_{x_{i_1+1/2}}^{x_{i_1+1/2}+\delta_{i_1+1/2,i_2}}R[V^n](x)dx - \right. \\
&& \left. - \int_{x_{i_2-1/2}}^{x_{i_2+1/2}}\int_{x_{i_1+1/2}-\delta_{i_1+1/2,i_2}}^{x_{i_1+1/2}}R[V^n](x)dx \right),\nonumber
\end{eqnarray}
and the displacement $\delta_{i_1+1/2,i_2}$ defined, at the centers of E and W facets of a square cell, by
\begin{equation*}
\delta_{i_1+1/2,i_2} = \sqrt{4\Delta t\>\nu_1(x_{i_1+1/2,i_2},t^n)}.
\end{equation*}
Note that the integrals appearing in all the fluxes of the form~\eqref{flux2d} have an integration 
domain consisting of a row (or column) of grid cells, among which at most one is intersected.
 
   \section{Numerical experiments.}
\label{tests}

   Several numerical experiments have been carried out with simple implementations of the FFSL method proposed above,
in order to assess its accuracy and stability features also in cases more complex  than those allowing a complete theoretical analysis. The accuracy of the proposed discretization has been evaluated against analytic solutions or reference solutions obtained by alternative discretizations in space and time.
Due to the close relationship between the methods
proposed here and those in~\cite{bonaventura:2014}, the choice of the test cases follows closely the outline in our previous paper,
in order to allow for a clear comparison between the conservative and nonconservative SL discretizations.

 \subsection{Constant coefficient case.}
\label{constcoeff}
In a first set of numerical experiments, the constant coefficient
 diffusion equation
\[
u_t   = \nu u_{xx} \ \ \ \ \ x\in[0,L]
\]
was considered, on an interval $[0,L] $ with $L=10.$  Periodic boundary conditions
were assumed and a Gaussian profile centered at $L/2$ was considered as the initial condition.
In this case, the exact solution  can be computed up to machine accuracy by separation
of variables and computation of its Fourier coefficients on a discrete mesh of
$N$ points with spacing $\Delta x=L/N.$ We consider the FFSL method
described in the previous Sections  on a time interval $[0,T] $ with $T=2, $
with time steps defined as $\Delta t=T/M.$ The   stability
parameter of standard explicit discretizations of the diffusion operator is defined
as   $\mu= \nu \Delta t/ \Delta x^2. $  We consider  the case
with   $\nu=0.05 $ first, whose relative errors  are reported in Tables~\ref{table1}-\ref{table2}, 
as computed in the $l_2$ and $l_{\infty}$, respectively.
The parallel results for the advection diffusion case 
\[
u_t  + au_x= \nu u_{xx} \ \ \ \ \ x\in[0,L]
\]
with $a=1.5, \nu=0.05$ are reported in Tables~\ref{table3}-\ref{table4}, respectively. 
In this case, the Courant number is defined as $ C=a\Delta t/ \Delta x $ and
the advective flux was computed using either linear or quadratic reconstruction,
along the lines of the well known method~\cite{lin:1996}. The fractional flux employed
yields in the low Courant number case a discretization that is equivalent to a flux
form Lax-Wendroff method.

   \begin{table}%[htbc]
 \centering
\begin{tabular}{||c|c|c|c|c|c|c||} 
\hline
\multicolumn{3}{||c|}{Resolution}  & \multicolumn{2}{c|}{$l_2$ rel. error (SL)} & \multicolumn{2}{c||}{$l_2$ rel. error (FFSL)} \\
\hline
$N$ & $M$ & $\mu$ &$I_1$ & $I_3$ & $R_0$ & $R_2$ \\
\hline 
200 & 50 & 0.8 &$1.50\cdot 10^{-2}$& $3.33\cdot 10^{-4}$ & $1.5\cdot 10^{-2} $& $3.79\cdot 10^{-4}$     \\
\hline 
200 & 100 & 0.4 &$1.45\cdot 10^{-2} $ & $2.91\cdot 10^{-4}$ & $1.45\cdot 10^{-2} $ & $2.02\cdot 10^{-4}$     \\
\hline
200 & 200 & 0.2  &$6.53\cdot 10^{-2}$ & $5.89\cdot 10^{-4}$& $6.53\cdot 10^{-2}$ & $2.28\cdot 10^{-4}$      \\
\hline
400 & 100 & 1.6 &$6.53\cdot 10^{-3}$ & $1.92\cdot 10^{-4}$ & $6.53\cdot 10^{-3}$ & $1.78\cdot 10^{-4}$    \\
\hline
400 & 200 & 0.8 &$1.48\cdot 10^{-2} $& $8.29\cdot 10^{-5}$ & $1.48\cdot 10^{-2}$ & $9.45\cdot 10^{-5}$    \\
\hline
400 & 400 & 0.4 &$1.43\cdot 10^{-2} $ & $7.43\cdot 10^{-5}$ & $1.43\cdot 10^{-2}$ & $5.03\cdot 10^{-5}$    \\
\hline
\end{tabular}
\caption{Relative errors for the  constant coefficient pure diffusion case in the $l_2$-norm, nonconservative (SL) and conservative (FFSL) scheme, first and third order space discretizations.}
 \label{table1}   
 \end{table}

   \begin{table}%[htbc]
 \centering
\begin{tabular}{||c|c|c|c|c|c|c||} 
\hline
\multicolumn{3}{||c|}{Resolution}  & \multicolumn{2}{c|}{$l_{\infty}$ rel. error (SL)} & \multicolumn{2}{c||}{$l_{\infty}$ rel. error (FFSL)} \\
\hline
$N$ & $M$ & $\mu$ &$I_1$ & $I_3$ & $R_0$ & $R_2$ \\
\hline 
200 & 50 & 0.8  &$1.73\cdot 10^{-2} $& $3.75\cdot 10^{-4}$ & $1.73\cdot 10^{-2} $& $4.44\cdot 10^{-4}$     \\
\hline 
200 & 100 &  0.4 &$1.67\cdot 10^{-2} $ & $6.46\cdot 10^{-4}$ & $1.67\cdot 10^{-2}$ & $2.36\cdot 10^{-4}$     \\
\hline
200 & 200 & 0.2 &$7.41\cdot 10^{-2} $& $1.63\cdot 10^{-3}$& $7.41\cdot 10^{-2}$& $2.67\cdot 10^{-4}$      \\
\hline
400 & 100 &  1.6 &$7.56\cdot 10^{-3}$ & $2.85\cdot 10^{-4}$ & $7.56\cdot 10^{-3}$ & $2.08\cdot 10^{-4}$    \\
\hline
400 & 200 & 0.8 &$1.71\cdot 10^{-2}$ & $8.86\cdot 10^{-5}$ & $1.70\cdot 10^{-2}$ & $1.11\cdot 10^{-4}$    \\
\hline
400 & 400 & 0.4 &$1.65\cdot 10^{-2} $ & $1.94\cdot 10^{-4}$ & $1.65\cdot 10^{-2}$ & $5.89\cdot 10^{-5}$    \\
\hline
\end{tabular}
\caption{Relative errors for the  constant coefficient pure diffusion case in the $l_{\infty}$-norm, nonconservative (SL) and conservative (FFSL) scheme, first and third order space discretizations.}
 \label{table2}   
 \end{table} 
 
    \begin{table}%[htbc]
 \centering
\begin{tabular}{||c|c|c|c|c|c|c|c||} 
\hline
\multicolumn{4}{||c|}{Resolution}  & \multicolumn{2}{c|}{$l_2$ rel. error (SL)} & \multicolumn{2}{c||}{$l_2$ rel. error (FFSL)} \\
\hline
$N$ & $M$ & C& $\mu$ &$I_1$ & $I_3$ & $R_0$ & $R_2$ \\
\hline 
200 & 50 & 1.2 & 0.8 &$1.21\cdot 10^{-2}$& $3.77\cdot 10^{-4}$ & $2.67\cdot 10^{-2} $& $6.46\cdot 10^{-4}$     \\
\hline 
200 & 100 & 0.6 &0.4 &$3.38\cdot 10^{-2} $ & $2.91\cdot 10^{-4}$ & $4.83\cdot 10^{-2} $ & $1.89\cdot 10^{-3}$     \\
\hline
200 & 200 & 0.3 &0.2  &$4.15\cdot 10^{-2}$ & $1.57\cdot 10^{-3}$& $1.11\cdot 10^{-1}$ & $2.67\cdot 10^{-3}$      \\
\hline
400 & 100 & 1.2 &1.6 &$5.02\cdot 10^{-3}$ & $1.84\cdot 10^{-4}$ & $1.25\cdot 10^{-2}$ & $2.14\cdot 10^{-4}$    \\
\hline
400 & 200 & 0.6 &0.8 &$1.30\cdot 10^{-2} $& $1.16\cdot 10^{-4}$ & $3.21\cdot 10^{-2}$ & $4.79\cdot 10^{-4}$    \\
\hline
400 & 400 & 0.3 & 0.4 &$2.95\cdot 10^{-2} $ & $4.27\cdot 10^{-4}$ & $4.41\cdot 10^{-2}$ & $6.67\cdot 10^{-4}$    \\
\hline
\end{tabular}
\caption{Relative errors for the  constant coefficient advection diffusion case in the $l_2$-norm, nonconservative (SL) and conservative (FFSL) scheme, first and third order space discretizations.}
 \label{table3}   
 \end{table} 
 
  \begin{table}%[htbc]
 \centering
\begin{tabular}{||c|c|c|c|c|c|c|c||} 
\hline
\multicolumn{4}{||c|}{Resolution}  & \multicolumn{2}{c|}{$l_{\infty}$ rel. error (SL)} & \multicolumn{2}{c||}{$l_{\infty}$ rel. error (FFSL)} \\
\hline
$N$ & $M$ & C& $\mu$ &$I_1$ & $I_3$ & $R_0$ & $R_2$ \\
\hline 
200 & 50 & 1.2 & 0.8 &$1.39\cdot 10^{-2}$& $4.85\cdot 10^{-4}$ & $3.05\cdot 10^{-2} $& $6.06\cdot 10^{-4}$     \\
\hline 
200 & 100 &0.6 &0.4 &$3.87\cdot 10^{-2} $ & $2.91\cdot 10^{-4}$ & $5.50\cdot 10^{-2} $ & $1.91\cdot 10^{-3}$     \\
\hline
200 & 200 & 0.3 &0.2  &$4.73\cdot 10^{-2}$ & $1.71\cdot 10^{-3}$& $1.28\cdot 10^{-1}$ & $2.69\cdot 10^{-3}$      \\
\hline
400 & 100 & 1.2&1.6 &$5.82\cdot 10^{-3}$ & $2.56\cdot 10^{-4}$ & $1.45\cdot 10^{-2}$ & $2.39\cdot 10^{-4}$    \\
\hline
400 & 200 & 0.6 &0.8 &$1.49\cdot 10^{-2} $& $2.51\cdot 10^{-4}$ & $3.68\cdot 10^{-2}$ & $4.85\cdot 10^{-4}$    \\
\hline
400 & 400 & 0.3 & 0.4 &$3.38\cdot 10^{-2} $ & $4.81\cdot 10^{-4}$ & $5.03\cdot 10^{-2}$ & $6.73\cdot 10^{-4}$    \\
\hline
\end{tabular}
\caption{Relative errors for the  constant coefficient advection diffusion case in the $l_{\infty}$-norm, nonconservative (SL) and conservative (FFSL) scheme, first and third order space discretizations.}
 \label{table4}   
 \end{table} 
 
 \noindent
 It can be observed that, in the linear case, the errors obtained by the SL and FFSL methods are in general of the same
 order of magnitude, with consistently smaller values for the conservative variant. The slightly 
larger errors of the FFSL case in the advection diffusion cases are partly explained by the fact that, in the present
implementation, for simplicity a second order reconstruction was employed for the advective terms, compared to
 the third order accurate reconstruction employed for the diffusion terms.

 \subsection{Linear diffusion with variable coefficients.}
\label{linvar}

\begin{figure}[htbc]
\begin{center}
\includegraphics[width=0.7\textwidth]{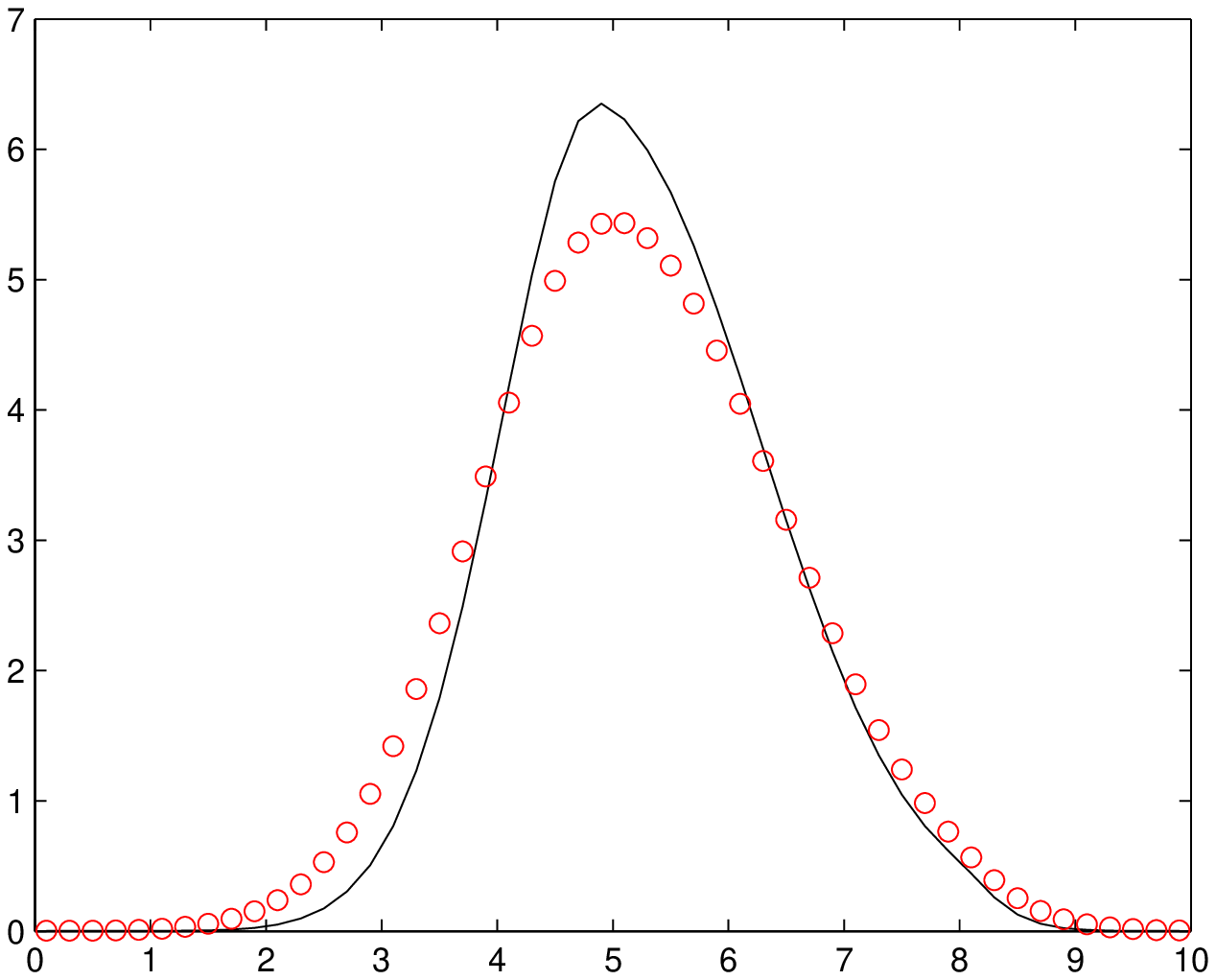}(a)
\includegraphics[width=0.7\textwidth]{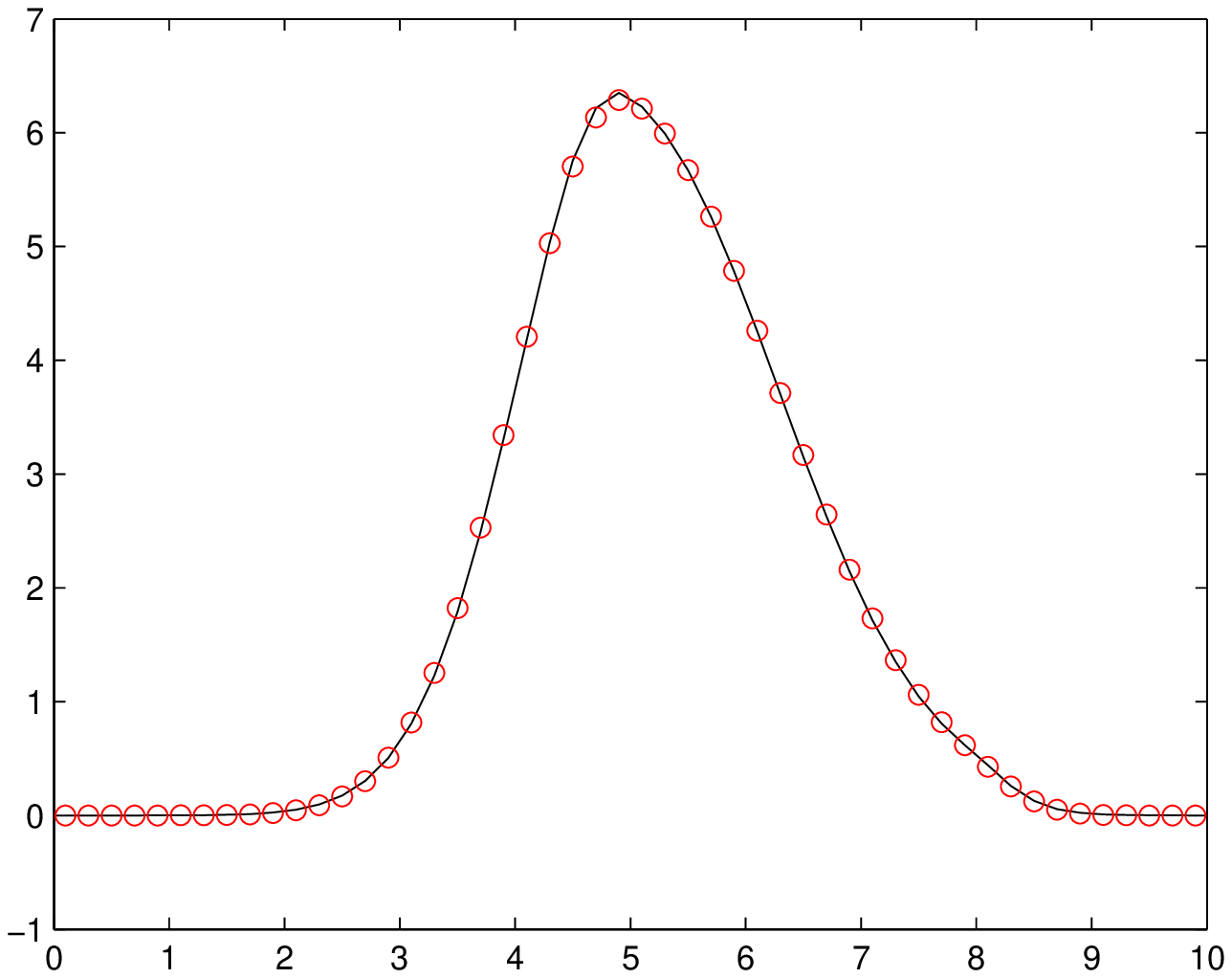}(b)
\caption{One-dimensional case with variable coefficients: numerical solution by FFSL method with (a) linear  (b) cubic reconstruction.
Numerical solution is represented by circles, reference solution by continuous line).}
\label{variable_ffsl}
\end{center}
\end{figure}

\begin{figure}[htbc]
\begin{center}
\includegraphics[width=0.7\textwidth]{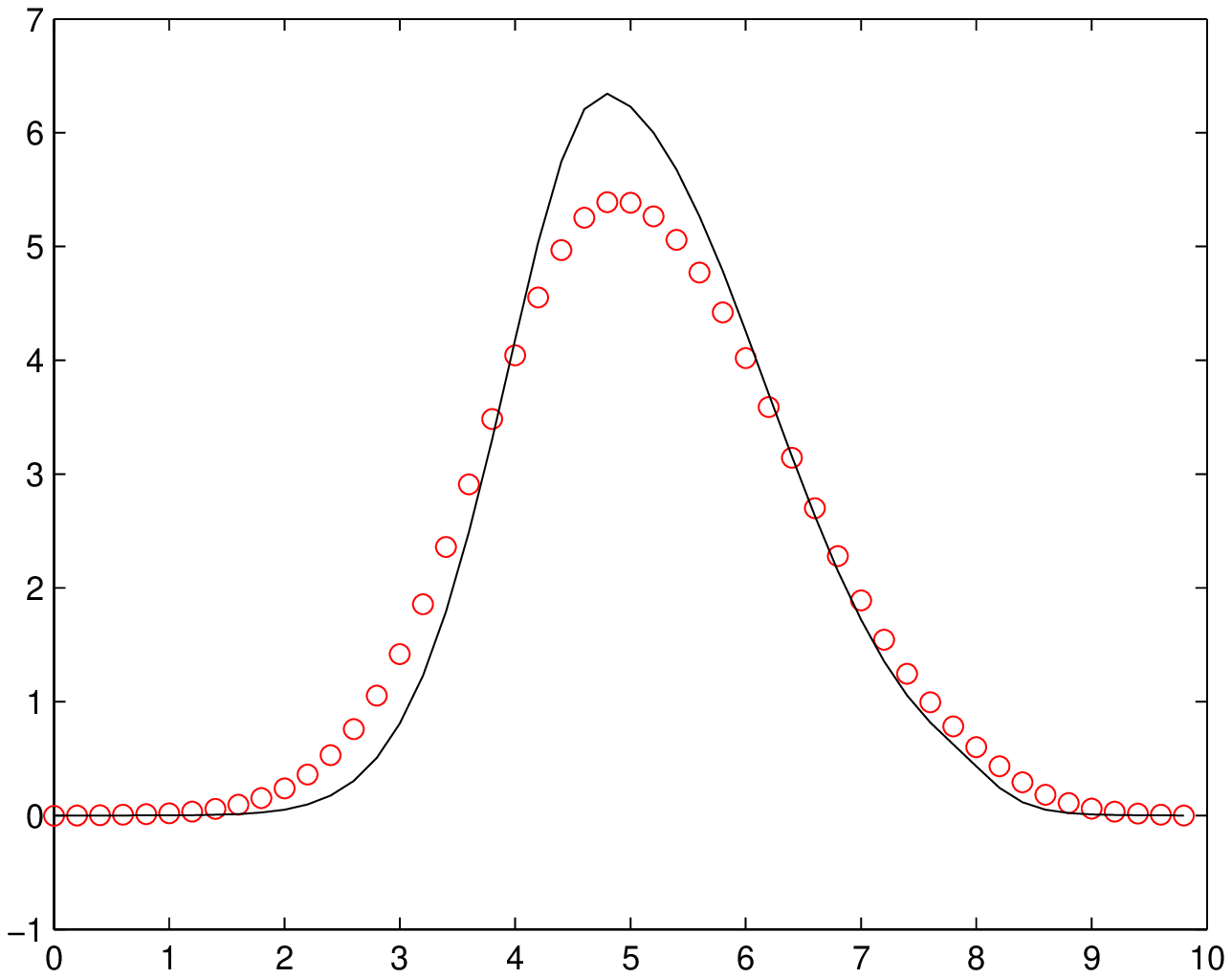}(a)
\includegraphics[width=0.7\textwidth]{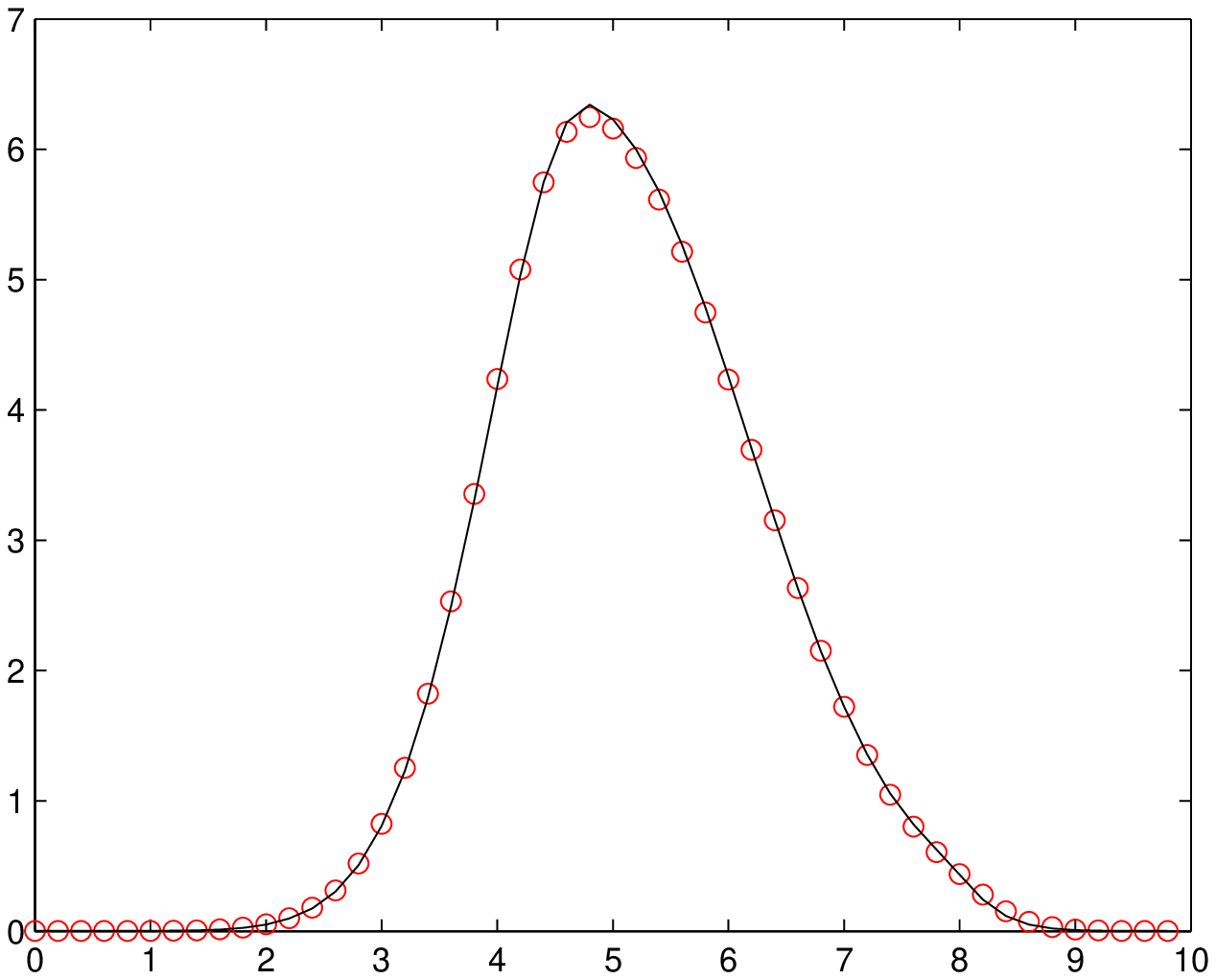}(b)
\caption{One-dimensional case with variable coefficients: numerical solution by SL method with (a) linear  (b) cubic reconstruction.
Numerical solution is represented by circles, reference solution by continuous line).}
\label{variable_sl}
\end{center}
\end{figure}

The  diffusion equation 
\[
u_t   = (\nu(x,t) u_{x} )_x\ \ \ \ \ x\in[0,L]
\]
was then considered  on a space interval $[0,L] $  and time interval $[0,T] $ with $L=10$  and  
 with $T=4,$
with time steps defined as $\Delta t=T/M.$ 
 Periodic boundary conditions were assumed and a Gaussian profile centered at $L/5$ was considered as the initial condition.
The  diffusivity field was given by

$$
\nu(x,t)= \frac 1{20}+ \frac 1{5} \xi(x)\sin\left(\frac{2\pi t}T\right)^2 ,
$$
respectively, where $ \xi(x) $ denotes the characteristic function of the interval $[0.5L,0.8L].$ This choice highlights the possibility to use the
proposed method seamlessly also with strongly varying diffusion coefficients.
In this case, no exact solution is available and reference solutions were computed using
the finite difference method described in the previous section with a four times higher spatial
resolution, coupled to a high order  multistep stiff solver
in which a small tolerance and maximum time step value were enforced. 
A plot of the numerical solutions obtained at the final time $T$ are displayed in Figures~\ref{variable_ffsl}-\ref{variable_sl}, for the
FFSL and SL method, respectively.

A more quantitative assessment of the FFSL solution accuracy can be gathered from Tables
~\ref{table5}-~\ref{table6}.
It can be observed that the FFSL method is always of slightly higher accuracy than the corresponding SL variant
at equivalent resolution. In all these computations, FFSL maintains mass conservation up to machine accuracy by construction, while
the average conservation error of the SL method is of the order of $10^{-3} $ of the total initial mass.

   \begin{table}%[htbc]
 \centering
\begin{tabular}{||c|c|c|c|c|c|c||} 
\hline
\multicolumn{3}{||c|}{Resolution}  & \multicolumn{2}{c|}{$l_2$ rel. error (SL)} & \multicolumn{2}{c||}{$l_2$ rel. error (FFSL)} \\
\hline
$N$ & $M$ & $\mu$ &$I_1$ & $I_3$ & $R_0$ & $R_2$ \\
\hline 
50 & 50 & 0.1 &$1.34\cdot 10^{-1}$& $1.10\cdot 10^{-3}$ & $1.30\cdot 10^{-1} $& $7.31\cdot 10^{-3}$     \\
\hline 
100 & 25 & 0.8 &$2.10\cdot 10^{-2} $ & $5.74\cdot 10^{-3}$ & $2.10\cdot 10^{-2} $ & $4.75\cdot 10^{-3}$     \\
\hline
100 & 100 & 0.2  &$7.26\cdot 10^{-2}$ & $4.21\cdot 10^{-3}$& $7.21\cdot 10^{-2}$ & $3.52\cdot 10^{-3}$      \\
\hline
200 & 50 & 1.6 &$1.02\cdot 10^{-2}$ & $4.68\cdot 10^{-3}$ & $9.57\cdot 10^{-3}$ & $3.00\cdot 10^{-3}$    \\
\hline
\end{tabular}
\caption{Relative errors for the  variable coefficient pure diffusion case in the $l_2$-norm, nonconservative (SL) and conservative (FFSL) scheme, first and third order space discretizations.}
 \label{table5}   
 \end{table}

   \begin{table}%[htbc]
 \centering
\begin{tabular}{||c|c|c|c|c|c|c||} 
\hline
\multicolumn{3}{||c|}{Resolution}  & \multicolumn{2}{c|}{$l_2$ rel. error (SL)} & \multicolumn{2}{c||}{$l_2$ rel. error (FFSL)} \\
\hline
$N$ & $M$ & $\mu$ &$I_1$ & $I_3$ & $R_0$ & $R_2$ \\
\hline 
50 & 50 & 0.1 &$1.50\cdot 10^{-1}$& $1.50\cdot 10^{-2}$ & $1.47\cdot 10^{-1} $& $1.28\cdot 10^{-2}$     \\
\hline 
100 & 25 & 0.8 &$2.72\cdot 10^{-2} $ & $7.97\cdot 10^{-3}$ & $2.77\cdot 10^{-2} $ & $8.21\cdot 10^{-3}$     \\
\hline
100 & 100 & 0.2  &$8.78\cdot 10^{-2}$ & $7.42\cdot 10^{-3}$& $8.87\cdot 10^{-2}$ & $6.71\cdot 10^{-3}$      \\
\hline
200 & 50 & 1.6 &$1.36\cdot 10^{-2}$ & $6.95\cdot 10^{-3}$ & $1.43\cdot 10^{-2}$ & $5.60\cdot 10^{-3}$    \\
\hline
\end{tabular}
\caption{Relative errors for the  variable coefficient pure diffusion case in the $l_{\infty}$-norm, nonconservative (SL) and conservative (FFSL) scheme, first and third order space discretizations.}
 \label{table6}   
 \end{table}

\subsection{Gas flow in porous media.}
\label{porous}
We reconsider the nonlinear  example, already treated in~\cite{bonaventura:2014}, of the one-dimensional equation of gases in porous media,
\begin{equation}\label{eq:mp}
u_t = \left(mu^{m-1}u_x\right)_x,
\end{equation}
focusing in particular on the so-called Barenblatt--Pattle self-similar solutions (see, e.g.,~\cite{barenblatt1952self}), which can be written in the form
\begin{equation}
u(x,t) = (t+t_0)^{-k} \left(A^2 - \frac{k(m-1)|x|^2}{2m(t+t_0)^{2k}}\right)_+^\frac{1}{m-1}
\end{equation}
where $t_0>0$, $A$ is an arbitrary nonzero constant and
$k = 1/(m+1).$
In order to adapt the FFSL scheme to this case, we recall what has been already observed for the nonconservative scheme in~\cite{bonaventura:2014}: since the solution (and hence, the diffusivity) may have a bounded support, a straightforward extension of~\eqref{delta_t} in the linearized form
\begin{equation}\label{delta_t_nl}
\delta^\pm_i = \sqrt{2\Delta t\>\nu\left(I[V^n](x_i\pm\delta^\pm_i)\right)}
\end{equation}
(where in our case $\nu(u)=mu^{m-1}$) would have multiple solutions for $x_i$ out of the support but in its neighbourhood -- to ensure a correct propagation of the solution, $\delta^\pm_i$ should be defined as the largest of such values. The same caution should be applied in extending the definition~\eqref{delta_ff}. A possible answer in this respect is to compute the $\delta^\pm_k$ via~\eqref{delta_t_nl}, then define
$\delta_k = (\delta^+_k+\delta^-_k)/2, $
where, of course, $x_k$ is now a cell interface.

\begin{figure}%[htbc]
\begin{center}
\includegraphics[width=0.65\textwidth]{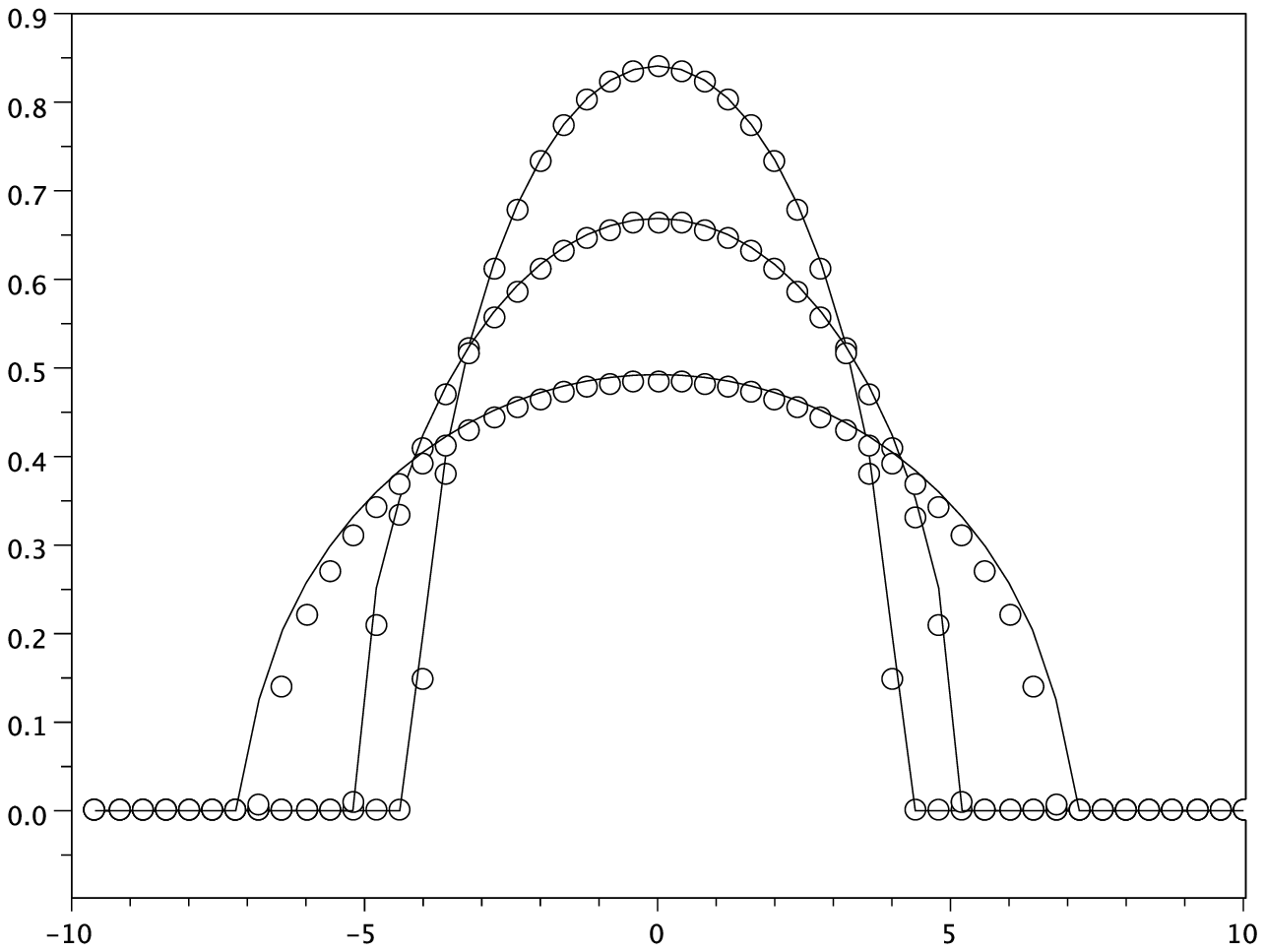}(a)
\includegraphics[width=0.65\textwidth]{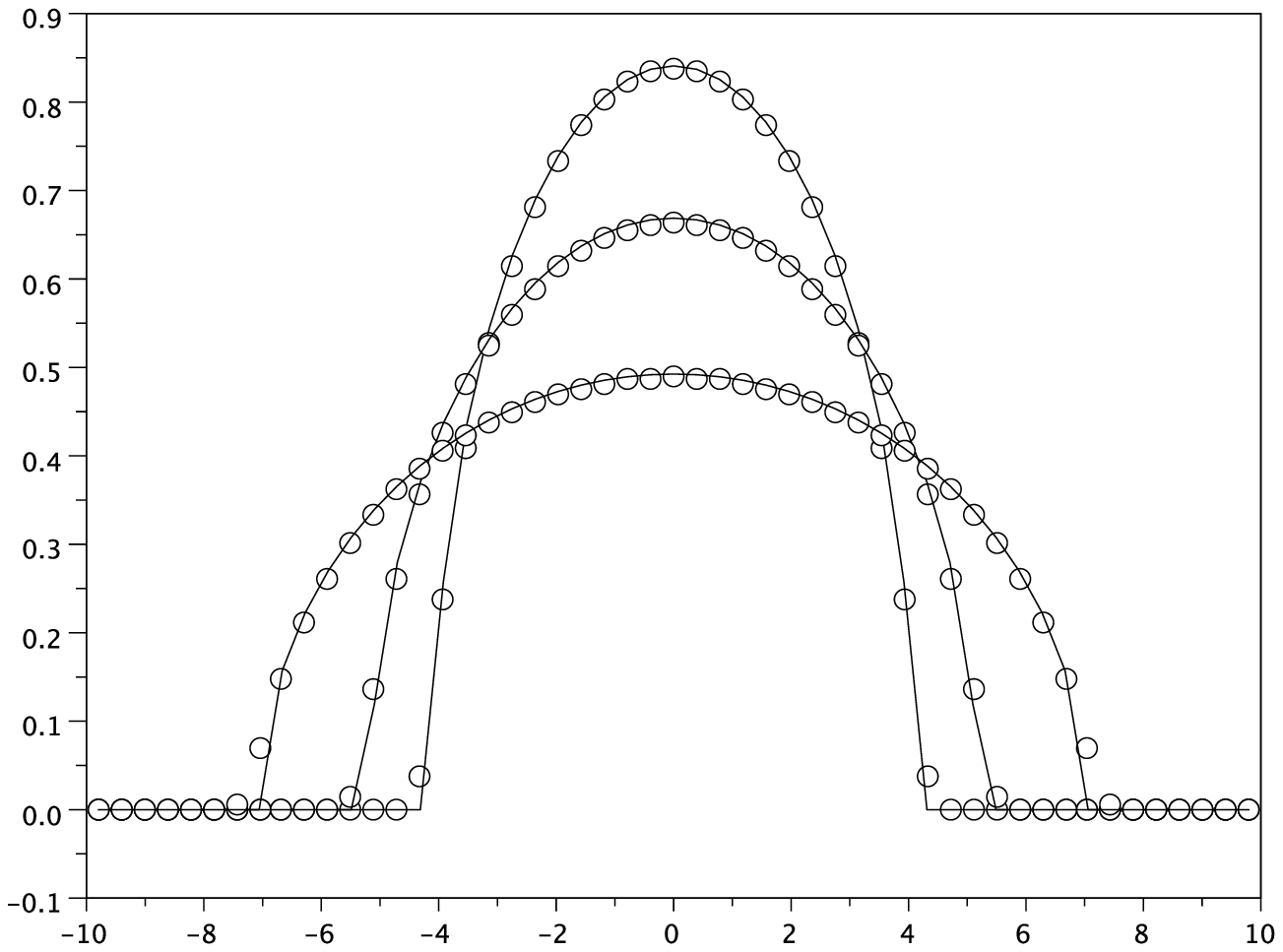}(b)
\caption{Evolution of the exact (continuous line) versus numerical (circles) Barenblatt--Pattle solution for T=1,4,16, nonconservative (a) and conservative (b) scheme.}
\label{fig:mp1}
\end{center}
\end{figure}

Figure~\ref{fig:mp1} compares the exact and approximate evolution of the Barenblatt--Pattle solution for~\eqref{eq:mp}, with $m=3$, $A=1$ and $t_0=1$. Approximate solutions have been computed with the nonconservative (a) and the conservative (b) scheme at $T=1,4,16$ on a mesh composed of 50 nodes with $\Delta t=0.05$, using respectively cubic interpolation and quadratic reconstruction. Note that the mass conservation constraint apparently improves the accuracy of the scheme, especially for larger simulation times.
This behaviour is confirmed, in terms of both absolute accuracy and convergence rate, by Table~\ref{table_bp}, which shows the numerical errors for the two schemes in the 2-norm, at $T=16$, under a linear refinement law. SL scheme has been tested with linear ($I_1$) and cubic ($I_3$) interpolation, whereas FFSL scheme with piecewise constant ($R_0$) and piecewise quadratic ($R_2$) reconstruction.

 \begin{table}%[htbc]
 \centering
\begin{tabular}{||c|c|c|c|c|c||} 
\hline
\multicolumn{2}{||c|}{Resolution} & \multicolumn{2}{c|}{$l_2$ relative error (SL)} & \multicolumn{2}{c||}{$l_2$ relative error (FFSL)} \\
\hline
$N$ & $M$ & $I_1$ & $I_3$ & $R_0$ & $R_2$ \\
\hline 
50 & 320 & 0.316 & $8.69\cdot 10^{-2}$ & 0.270 & $2.67\cdot 10^{-2}$     \\
\hline
100 & 640 &  0.212 & $4.76\cdot 10^{-2}$& 0.156 & $1.65\cdot 10^{-2}$      \\
\hline
200 & 1280 &  0.171 & $4.84\cdot 10^{-2}$ & $9.85\cdot 10^{-2}$ & $1.09\cdot 10^{-2}$    \\
\hline
400 & 2560 & 0.135 & $4.49\cdot 10^{-2}$ & $6.03\cdot 10^{-2}$ & $2.69\cdot 10^{-3}$    \\
\hline
800 & 5120 & 0.107 & $3.22\cdot 10^{-2}$ & $3.72\cdot 10^{-2}$ & $2.79\cdot 10^{-3}$   \\
\hline
\end{tabular}
\caption{Relative errors for the Barenblatt--Pattle solution in the 2-norm, nonconservative (SL) and conservative (FFSL) scheme, first and third order space discretizations.}
 \label{table_bp}   
 \end{table}

\subsection{Variable coefficient case in two space dimensions.}
\label{varcoeff2}

\begin{figure}%[htbc]
\begin{center}
\includegraphics[width=0.75\textwidth]{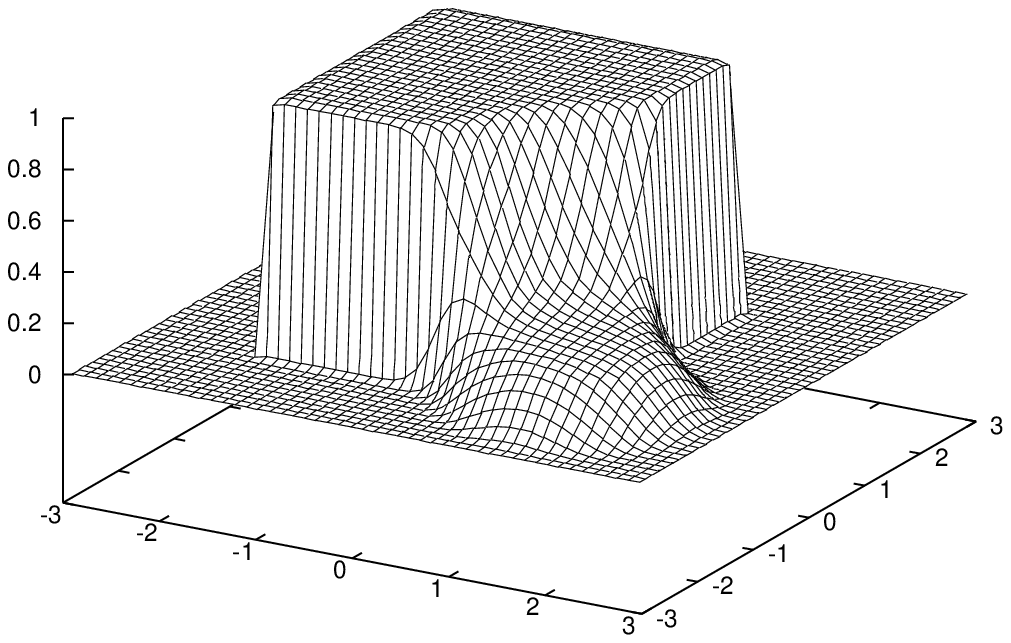}(a)
\includegraphics[width=0.65\textwidth]{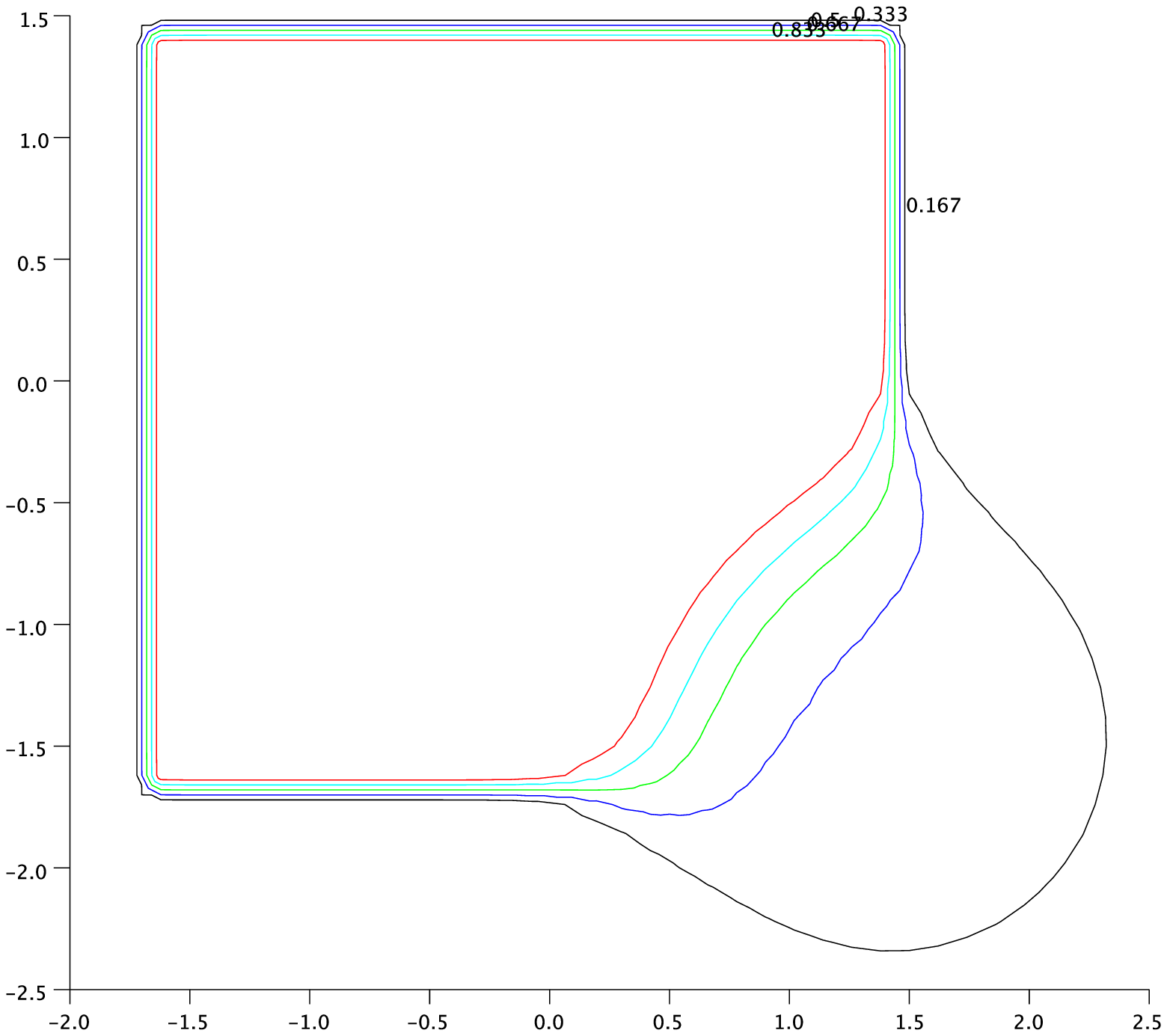}(b)
\caption{Variable isotropic diffusion, graph  (a) and level curves (b) of the solutions.}
\label{fig:2d}
\end{center}
\end{figure}

As a two-dimensional  test, we consider the equation
\[
u_t = \text{div}(\nu(x)\nabla u)
\]
on $\Omega=[-3,3]^2$ with periodic boundary conditions. 
The initial condition given by the characteristic function of the square $\Sigma=[-1.5,1.5]^2$. The isotropic diffusivity $\nu$ is given by
\[
\nu(x) = e^{-5|x-x_0|^2},
\]
where $x_0=(1.5,-1.5).$ The diffusion is therefore concentrated in a corner on the boundary of the set $\Sigma$. The effect of this diffusion is to move mass from the interior of $\Sigma$ to the exterior, in the neighbourhood of the point $x_0$. Figure~\ref{fig:2d} shows the numerical solution along with its level curves at $T=2$, with a $50\times 50$ space grid, $R_0$ reconstruction and time step $\Delta t =0.05$. Note that, despite being conceptually obtained by directional splitting, the two-dimensional scheme ~\eqref{eq:schema2d} does not suffer from anisotropies induced by grid orientation.

 \section{Conclusions.}
\label{conclu}
We have extended the SL approach of ~\cite{bonaventura:2014} to achieve a  fully conservative, flux form
discretization for linear and nonlinear parabolic problems. 
A consistency and stability analysis of the method has also been presented, along with a strategy
to couple the FFSL discretizations of advection and diffusion terms. 
A number of numerical simulations validate the proposed method,
showing that it is equivalent in accuracy to its nonconservative variant, while allowing to
maintain mass conservation at machine accuracy. The proposed method
 could represent an important complement to the many conservative, flux-form semi-Lagrangian schemes for advection
 proposed in the literature. Future developments will focus on the
 application of the proposed approach to  conservation laws with nonlinear parabolic terms, such as e.g. the Richards
 equation, and to the extension of this technique to high order discontinuous finite elements
 discretizations such as those proposed in ~\cite{restelli:2006}, ~\cite{tumolo:2013}.

\section*{Acknowledgements.}

This research work has been financially supported by the INDAM--GNCS project
\textit{Metodi numerici semi-impliciti e semi-Lagrangiani per sistemi iperbolici di leggi di bilancio}  
and by Politecnico di Milano.

\bibliographystyle{plain}
\bibliography{diff_ffsl}

\end{document}